\documentclass[a4paper,10pt]{article}%

\usepackage[utf8]{inputenc}%
\usepackage[T1]{fontenc}%
\usepackage[english]{babel}%

\usepackage{lmodern}%
\usepackage[final]{graphicx}%
\usepackage{amsmath,amssymb,amsthm}%
\usepackage{hyperref}%
\usepackage{enumitem}%
\usepackage[labelfont=bf,labelsep=period]{caption}%

\allowdisplaybreaks%
\newcommand*{\Rbb}{\mathbb{R}}% reals
\newcommand*{\Zbb}{\mathbb{Z}}% integers

\newcommand*{\Gc}{\mathcal{G}}%

\newcommand*{\etab}{{\boldsymbol\eta}}%
\newcommand*{\Gmb}{{\boldsymbol\Gamma}}%
\newcommand*{\vphi}{\varphi}%
\newcommand*{\pib}{{\boldsymbol\pi}}%
\newcommand*{\om}{\omega}%
\newcommand*{\omb}{{\boldsymbol\omega}}%
\newcommand*{\Om}{\Omega}%

\newcommand*{\ul}{u_{-}}% left density
\newcommand*{\ur}{u_{+}}% right density

\newcommand*{\emp}[1]{\emph{#1}}%

\numberwithin{equation}{section}%

\theoremstyle{plain}%
\newtheorem{theorem}{Theorem}%

\theoremstyle{definition}%

\DeclareMathOperator{\Exp}{\mathbf{E}}%
\DeclareMathOperator{\infl}{\mathrm{infl}}%

\makeatletter
\let\@fnsymbol\@alph%\@arabic
\makeatother

\newcommand{\BM}{\url{m.balazs@bristol.ac.uk}, School of Mathematics, University of Bristol, UK.\ }%
\newcommand{\NAL}{\url{attilalaszlo.nagy@gmail.com} or \url{nagyal@math.bme.hu}, Institute of Mathematics, Budapest University of Technology and Economics, HU.\ }%
\newcommand{\TB}{\url{balint.toth@bristol.ac.uk}, School of Mathematics, University of Bristol, UK. MTA-BME Stochastics Research Group and R\'enyi Institute Budapest, HU.\ }%
\newcommand{\TI}{\url{tothi@math.bme.hu}, Institute of Mathematics, Budapest University of Technology and Economics, HU.\ }%

\newcommand{\titl}{
Coexistence of shocks and rarefaction fans: complex phase diagram of a simple hyperbolic particle system
}%

\newcommand{\auths}{
M\'arton Bal\'azs\,\footnote{\BM}\\[1em]
\and
Attila L\'aszl\'o Nagy\,\footnote{\NAL}\\[1em]
\and
B\'alint  T\'oth\,\footnote{\TB}\\[1em]
\and
Istv\'an T\'oth\,\footnote{\TI}\\[1em]
}%

\newcommand{\dat}{\today}%

\title{\titl}%
\author{\auths}%
\date{\dat}%

\begin{document}
\maketitle%

\begin{abstract}
This paper investigates the non-equilibrium hydrodynamic behavior of a simple totally asymmetric interacting particle system of particles, antiparticles and holes on $\mathbb{Z}$. Rigorous hydrodynamic results apply to our model with a hydrodynamic flux that is exactly calculated and shown to change convexity in some region of the model parameters. We then characterize the entropy solutions of the hydrodynamic equation with step initial condition in this scenario which include various mixtures of rarefaction fans and shock waves. We highlight how the phase diagram of the model changes by varying the model parameters.
\end{abstract}

\bigskip


\bigskip

\noindent \textbf{Acknowledgement.}
M.\ B.\ and B.\ T.\ acknowledge support from the grants OTKA K100473 and OTKA K109684. B.\ T.\ also thanks support from the Leverhulme International Network Grant ``Laplacians, Random Walks and Quantum Spin Systems''. We are indebted to an anonymous referee whose detailed and helpful comments greatly improved the exposition of the paper.

\section{Introduction}%

This paper focuses on a particular nearest-neighbor interacting particle system over $\Zbb$ with one conserved quantity. The model consists of particles, antiparticles and holes with the following types of moves:
\begin{enumerate}[label={(\arabic*)}]
\item \textbf{Exclusion.} (Anti)particles execute totally asymmetric exclusions to the (negative) positive direction, respectively.
\item \textbf{Annihilation.} An adjacent particle-antiparticle pair can annihilate each other producing two holes.
\item \textbf{Pair creation.} An antiparticle-particle pair can be created from two adjacent holes.
\end{enumerate}
The evolution is a continuous time Markov jump process the jump rates of which are chosen such that the model is \emp{attractive} and possesses \emp{product-form stationary distributions}. These choices greatly facilitate the analysis of the macroscopic behavior of the system. By rescaling time and space with the same factor (Eulerian scaling), we arrive to a hydrodynamic limit described by a nonlinear partial differential equation. As usual in this area nonlinearity comes from the \emp{hydrodynamic flux} which is roughly speaking the net flux of (anti)particles across a bond in equilibrium.

\bigskip
\noindent\emp{Earlier results and result of this paper.} In many well-known examples of particle systems, such as simple exclusion or constant rate zero range processes, the hydrodynamic flux function is proved to be strictly concave, sometimes strictly convex. Entropy solutions of the hydrodynamic equation in these cases are well known: either shock wave or rarefaction fan emerges from step initial data (also called Riemann problem). The first realistic model that exhibits non-convex nor-concave hydrodynamic flux appeared in the seminal paper by Katz, Lebowitz and Spohn \cite{kls} which in its one-dimensional form can be considered as a generalization of simple exclusion. The KLS model, though it is non-attractive, has become popular and it has been under extensive studies during the past decades, see \cite{popkov,hager,chowd,popkovetal,junction}, showing spectacular properties in, e.g.,\ its phase diagram.

We highlight another generalization of simple exclusion, \emp{PushASEP} \cite{pushasep}, which produces non-convex nor-concave hydrodynamic flux. The nature of this model is slightly different from ours in that it requires non-nearest neighbor interactions and jump rates that depend on the configuration over the whole interval of jumps.

To the best of our knowledge the hydrodynamic behavior beyond the shock is not rigorously established for either of the known examples with non-convex nor-concave flux. In contrary, the model we introduce here is fully covered under the hydrodynamic results of the literature. It is an \emp{attractive} hyperbolic nearest neighbor particle system of just three possible states on a site which belongs to a more general family of misanthrope processes \cite{coco}. One of the main conclusions is that, while hydrodynamics beyond the shock is rigorously established for this very simple attractive model, it also exhibits a non-trivial flux of both concave and convex pieces in some parameter range. We have combined recent results of PDE solution theory (see \cite[Ch. 2]{holdenrisebro} and \cite{hayes,fossati}) with that of hydrodynamic limits established for general interacting particle systems \cite{bagurasa,bagurasastrong} to obtain rigorous results for the phase diagram of our model. We characterize all the possible cases of the Riemann problem which can consist of coexisting shock wave and rarefaction fan regions in various orders depending on the density values we started off from. We also show how the whole phase diagram evolves as one varies the underlying parameter values of the system. At some cumbersome calculations we used the aid of the computer. Finally, we mention that most of the results of the present paper originate from \cite{riemann}. As a nice add-on, we have developed some stochastic simulation programs (see \cite{sim}) for the evolution of our particle system.

\bigskip
\noindent\emp{Organization of the paper.} We define our microscopic model in Section \ref{sec:micromodel}. We calculate some of its key properties in Section \ref{sec:props}. The detailed investigation of the Riemann problem is contained in Section \ref{sec:entropysols}.

\section{The microscopic model}\label{sec:micromodel}%

The Markov process we consider consists of particles, antiparticles and holes interacting with each other on the one-dimensional integer lattice $\Zbb$. At lattice point $i\in\Zbb$, $\om_i$ will denote the presence of a particle ($\om_i=1$), that of an antiparticle ($\om_i=-1$) or the absence of them ($\om_i=0$). Thus our configuration space is $\Om:=\{-1,0,+1\}^{\Zbb}$ an element of which is denoted by  $\omb=(\om_i)_{i\in \Zbb}$. The continuous time Markovian jump dynamics we attach on top of $\Om$ is then made up of the moves
\[
\omb\;\longrightarrow\;\omb-\delta_j+\delta_{j+1}\in\Om,
\]
which happens at rate $p(\om_j,\om_{j+1})$ for every $j\in\Zbb$, where $\delta$ denotes the Kronecker symbol, i.e.\ $\delta_j(i)$ is $1$ if $i=j$ and zero otherwise. These changes take place independently conditioned on a given configuration $\omb\in\Om$. This process can easily be constructed in the usual way (see \cite{liggettbible}) and has infinitesimal generator $\Gc$ acting on a cylinder function (i.e., one that only depends on a finite number of coordinates) $\vphi:\Om\to\Rbb$ as
\begin{equation}\label{eq:infgen}
\big(\Gc\,\vphi\big)(\omb)=
\sum_{j\in\Zbb}\,p(\om_j,\om_{j+1})\cdot\big(\vphi(\omb-\delta_j+\delta_{j+1})-\vphi(\omb)\big)\qquad (\omb\in\Om).
\end{equation}
With a slight abuse of notation the state of the process at time $t\geq 0$ is also denoted by $\omb(t)=\big(\om_i(t)\big)_{i\in\Zbb}$. We now specify $p$ to be of the following special form. For $x,y\in\{-1,0,+1\}$ let
\begin{equation}\label{eq:modelrates}
p(x,y)=\left\{
\begin{array}{ll}
c, & \quad\hbox{$x=0$, $y=0$;} \\
a, & \quad\hbox{$x=+1$, $y=-1$;} \\
a\cdot \frac{1 - d}{2}, & \quad\hbox{$x=0$, $y=-1$;} \\
a\cdot \frac{1 + d}{2}, & \quad\hbox{$x=+1$, $y=0$;} \\
0, & \quad\hbox{otherwise,}
\end{array}
\right.
\end{equation}
where $a,c>0$ are (distinct) model parameters denoting the rates of \emp{annihilation} and \emp{pair creation}, respectively. With parameter $d\in(-1,1)$ we can adjust the symmetry of (anti)particles' jumping rate. In plain words, the above dynamics boils down to the following simple rules: two adjacent holes can produce an antiparticle-particle pair (in that order) with rate $c$, antiparticles and particles can only hop to the negative and positive directions with rates $a\cdot \frac{1-d}{2}$ and $a\cdot\frac{1+d}{2}$, respectively, and when a particle meets an antiparticle they annihilate each other with rate $a$. All other moves are suppressed.

By rescaling time by $a$, without loss of generality, we can assume that the annihilation rate is set to be $1$. Next, we invoke the definition of \emp{attractivity} from \cite[pp.\ 71--72]{liggettbible}. Formally speaking we require the relation $\Exp f(\omb(t))\leq \Exp f(\etab(t))$ to hold for all monotone functions $f:\Om\to\Rbb$ and for all $t\geq0$ whenever the initial configuration $\etab(0)$ of $\etab$ dominates that of $\omb$ in the coordinate-wise order. Along the lines of \cite[III., Thm.\ 2.2]{liggettbible}, attractivity is equivalent to $p$ being monotone non-decreasing (non-increasing) in its first (second) variable, respectively. This means that our dynamics is attractive if and only if $c\leq \frac{1-|d|}{2}$ which we assume throughout the paper. (We do not expect different behavior in the non-attractive regime either but hydrodynamics beyond the shock is not rigorously established there and the interesting phenomena of this paper happens in the attractive parameter domain anyway.)

The model outlined above is the simplest member, after simple exclusion, of the family of \emp{misanthrope processes} \cite{coco}. In the above particular form, it first appeared in \cite[pp.\ 179--180]{tothvalkoperturb}.

\section{Properties of the model}\label{sec:props}%

First, we introduce the one-parameter family of product-form extremal stationary distributions of the above model. Then in Subsections \ref{sec:hydro} and \ref{sec:flux} we establish some hydrodynamic properties of the model. These are then used in Section \ref{sec:entropysols} to determine the entropy solution of the Riemann problem.

Due to attractivity a very convenient feature of our model with rates \eqref{eq:modelrates} is that it possesses \emp{product-form stationary distributions}. We define
\begin{equation}\label{eq:paramb}
b :\,= \frac{1}{2 + \frac{1}{\sqrt{c}}},
\end{equation}
which is a bijection between the parameters $c\in(0,\frac{1}{2}]$ and $b\in(0,\frac{1}{2+\sqrt{2}}]$. Now, let $\theta\in\Rbb$ be a generic real that is also referred to as the \emp{chemical potential}. Then define the one-site marginal measure by
\begin{align}\label{eq:statmarginal}
\Gamma^{\theta} :\,= \big(\Gamma^{\theta}(-1),\, \Gamma^{\theta}(0),\, \Gamma^{\theta}(1)\big) = \bigg( \frac{b\, \exp(-\theta)}{Z(\theta)},\,\frac{1 - 2b}{Z(\theta)},\,\frac{b\, \exp(\theta)}{Z(\theta)} \bigg),
\end{align}
with \emp{partition sum} $Z(\theta) :\,= 1 + 2b\cdot(\cosh(\theta) - 1)$. The product measure $\Gmb^{\theta}$ on $\Om$ is built from these marginals: $\Gmb^{\theta} :\,= \bigotimes_{j\in \Zbb}\Gamma^{\theta}$. We denote the expectation with respect to this measure by $\Exp^{\theta}$.
The stationarity of the above measure was carried out in \cite{coco}. For the ergodicity part we refer to \cite[pp. 1350--1352, Sec. 3 and further references therein]{bagurasa}.

\subsection{Hydrodynamics}\label{sec:hydro}%

Eulerian hydrodynamic limits are well established for a large class of \emp{attractive} asymmetric particle systems \cite{rezakhanlou,landim,kipnislandim} or more generally \cite{bagurasa, bagurasastrong} that also includes our model. Below we recall some elements of this theory following the latter articles. Let $N$ be the rescaling factor and define
\[
\alpha^N(\mathrm{d}x, t) :\,= \frac{1}{N}\sum_{j\in\Zbb}\om_{j}(t\cdot N)\cdot\delta_{j/N}(\mathrm{d}x)
\]
as the \emp{rescaled empirical measure} of $\omb(t)$ at macroscopic space and time $x\in\Rbb$ and $t\geq0$, respectively.

Now, the sequence of initial probability distributions $(\pib_N)_{N\in\Zbb^+}$ can be chosen arbitrarily such that the empirical measure $\alpha^N$ of $\omb^N(0)$, distributed as $\pib_N$, converges in probability to $u_0(\,\cdot\,)\,\mathrm{d}x$, where $u_0$ is some deterministic bounded measurable profile on $\Rbb$ (see \cite[pp. 1346--1348, Sec. 2.3]{bagurasa}). Then for each $t>0$: the $\alpha^N$ of $\omb^N(t)$ converges to $u(\,\cdot\,,t)\,\mathrm{d}x$ in probability, where the density profile $u(\,\cdot\,,t):\Rbb\to[-1,1]$ is one of the weak solutions of the \emp{conservation law}
\begin{equation}\label{eq:hydrogeneral}
\begin{aligned}
\partial_t u + \partial_x G(u) &= 0;\\
u(\,\cdot\,,0) &= u_0(\cdot).
\end{aligned}
\end{equation}
In this PDE, the \emp{hydrodynamic flux} $G:[-1,1]\to\Rbb_0^+$ is
\[
G(\varrho) = \Exp_{\nu_{\varrho}} p(\om_0,\om_1)\qquad (\varrho\in[-1,1]),
\]
where the measure $\nu_{\varrho}$ is a \emp{translation-invariant extremal stationary distribution} of the process $(\omb(t))_{t\geq0}$ corresponding to density $\varrho$. In our case, as we discussed in the beginning of Section \ref{sec:props}, these measures can be expressed in product form. In plain words $G$ explains the net flux across a bond in equilibrium (recall the dynamics of $\omb$ being totally asymmetric) and it will be determined in Section \ref{sec:flux}.

Subsequently, we will only consider the \emp{Riemann problem} (step initial datum) in view of the initial value problem \ref{eq:hydrogeneral}, that is
\begin{equation}\label{eq:stepinitcond}
u_0(x) = \left\{
\begin{array}{ll}
\ul \quad &\mbox{if} \quad x \leq 0; \\
\ur \quad &\mbox{if} \quad x > 0, \\
\end{array}
\right.
\end{equation}
where $\ul\neq\ur\in[-1,1]$ are the (initial) densities on the left and the right-hand side of the origin. In case of the microscopic process for sake of simplicity we set the density values ($\ul,\ur$) by picking the appropriate site-dependent parameters $\theta_j$ of the product measure $\Gmb^{\theta}$:
\begin{equation*}
\Exp^{\theta_j}(\om_j) =
\left\{
\begin{array}{ll}
\ul \quad &\mbox{if} \quad j \leq 0;  \\
\ur \quad &\mbox{if} \quad j > 0. \\
\end{array}
\right.
\end{equation*}
As we discussed above the initial measure can be chosen from a more general set of measures. Indeed this choice of initial measures will not play any role in the present article.

It is more convenient to reparametrize the marginal $\Gamma^{\theta}$ by the (\emp{signed}) \emp{density} of particles instead of the chemical potential $\theta$. Let $v(\theta)$ be the expected number of particles occupying an arbitrary lattice point under the measure $\Gmb^{\theta}$:
\[
v(\theta) :\,= \Exp^{\theta}\om_0 = \frac{2b\,\sinh(\theta)}{1 + 2b\cdot(\cosh(\theta) - 1)} \qquad (\theta\in\Rbb).
\]
It is not hard to see that the assignment $\theta\mapsto v(\theta)\in[-1,1]$ is strictly monotone hence bijective. Its inverse $v\mapsto\theta(v)$ can explicitly be calculated:
\begin{equation}\label{eq:densityparam}
\theta(v)=
\log\bigg(\frac{(1 - 2b)\cdot v + \sqrt{4b^2 + (1-4b)\cdot v^2}}{2b\cdot(1-v)}\bigg) \qquad (v\in[-1,1]).
\end{equation}
With a slight abuse of notation we will freely switch between these two parametrizations of the stationary measure.

To close this section, we mention a well-known strategy for solving \eqref{eq:hydrogeneral}. When $G$ has continuous derivative and $u$ is smooth, one can rewrite \eqref{eq:hydrogeneral} as $\partial_t u + G'(u)\cdot \partial_x u = 0$. The \emp{characteristic curves} are then defined as lines of the $x$-$t$ space along which the solution is constant. To find these curves one solves $\frac{\mathrm{d}}{\mathrm{d}t}u(x(t),t)=0$ which translates to $\frac{\mathrm{d}x(t)}{\mathrm{d}t}=G'(u)$ after comparison to the original equation. The curve then can be traced back to the initial condition and takes the form $x = x_0 + t\cdot G'(u_0(x_0))$. The quantity $G'(u_0(x_0)))$ is known as the \emp{characteristic velocity}. Finally, to look up the value of the solution $u$ at an arbitrary point $(x,t)\in\Rbb\times\Rbb^+_0$ one just follows the characteristic curve back to the initial time and reads the value of $u_0$ there.

The previous program relies on $u$ being smooth which sometimes just fails to happen. That is there can exist points of $\Rbb\times\Rbb^+$ which more than one characteristic lines hit carrying different initial values. In these cases classical (differentiable) solutions of \eqref{eq:hydrogeneral} cease to exist and more general, so-called \emp{weak solutions}, need to be defined. By a careful selection of the weak solutions one can uniquely identify the physically relevant one which is also referred to as the \emp{entropy solution}. For more details and properties of the entropy solutions we refer to the monographs \cite{smoller,holdenrisebro}.

In our case the initial condition $u_0$ is the step function \eqref{eq:stepinitcond} hence the characteristic velocities can only attain two values: $G'(\ul)$ and $G'(\ur)$. Depending on the relation between the densities $\ul,\ur$ and on the monotonicity of $G'$ in $[\min(\ul,\ur),\max(\ul,\ur)]$ several possibilities can emerge. In the simplest case when $G$ is strictly concave one can essentially distinguish two cases:
(1) if $G'(\ul) > G'(\ur)$ then the characteristic lines meet in some finite time hence a \emp{shock wave}, that is a moving discontinuity appears in the entropy solution;
(2) if $G'(\ul) < G'(\ur)$ then the characteristic lines are moving away from each other giving rise to a \emp{rarefaction fan} in which the initial sharp discontinuity smears out in time.

In the general case, based on \cite[pp. 30--34 of Sec. 2.2]{holdenrisebro}, the unique entropy solution of \eqref{eq:hydrogeneral} with initial condition \eqref{eq:stepinitcond} can explicitly be given as
\begin{equation}\label{eq:generalentropysol}
u(x,t) =
\left\{
\begin{aligned}
&\ul &&
\text{if $x\leq (G_{\frown})'(\ul)\cdot t$};\\
&\big[(G_{\frown})'\big]^{-1}(x/t) &&
\text{if $(G_{\frown})'(\ul)\cdot t < x\leq (G_{\frown})'(\ur)\cdot t$};\\[0.3em]
&\ur &&
\text{if $x > (G_{\frown})'(\ur)\cdot t$},
\end{aligned}
\right.
\end{equation}
where $\ul>\ur$ and $G_{\frown}$ is the \emp{upper concave envelope} of $G$ defined to be the smallest concave function which is greater than or equal to $G$ above the interval $[\ur,\ul]$. The derivative of $G_{\frown}$ is denoted by $(G_{\frown})'$, while its inverse by $[(G_{\frown})']^{-1}$. Analogously, in the opposite case when $\ur>\ul$ the solution can be obtained by the help of the \emp{lower convex envelope} of $G$.

\subsection{Flux}\label{sec:flux}%

We explicitly calculate the \emp{hydrodynamic flux} $G$ of \eqref{eq:hydrogeneral} corresponding to our microscopic model. As we noted in the previous Subsection \ref{sec:hydro} the hydrodynamic flux is the expected current rate across a bond under the translation-invariant extremal stationary distribution $\Gmb^{\theta}$ (see the beginning of Section \ref{sec:props}). That is
\begin{align*}
G(v) = \Exp^{\theta(v)} p(\om_0,\om_1)
& = \frac{1 - d}{2}\cdot \Gamma^{\theta(v)}(0)\cdot\Gamma^{\theta(v)}(-1)
+ \frac{1 + d}{2}\cdot \Gamma^{\theta(v)}(1)\cdot\Gamma^{\theta(v)}(0)\\
& + c \cdot \Gamma^{\theta(v)}(0)\cdot\Gamma^{\theta(v)}(0)
+ 1 \cdot \Gamma^{\theta(v)}(1)\cdot\Gamma^{\theta(v)}(-1),
\end{align*}
where $v\in[-1,1]$. We remark that even the last term of the annihilation step counts as one towards the flux as only a unit signed charge is transferred over the bond in question.

\begin{figure}[!htb]
\centering
\includegraphics[width=0.75\textwidth]{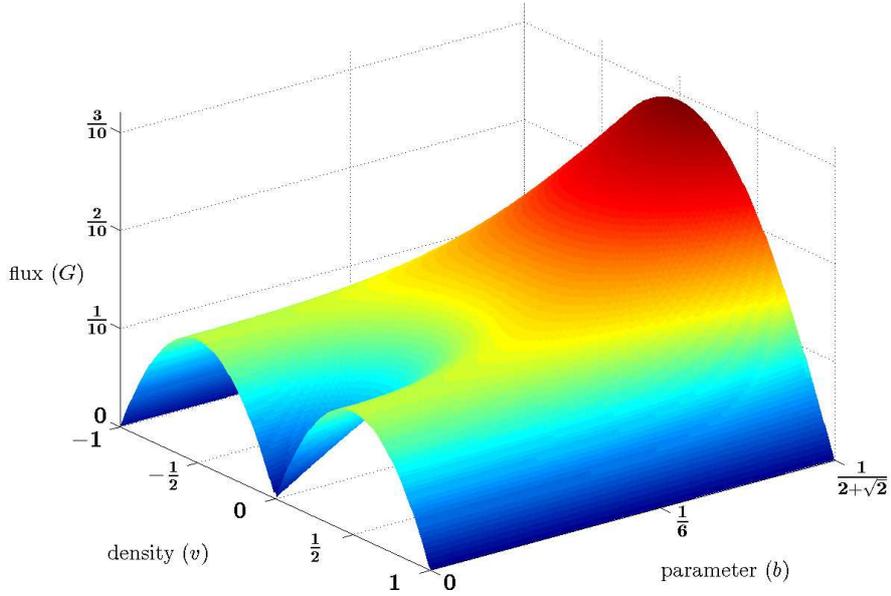}
\caption[Hydrodynamic flux]{The hydrodynamic flux $G$ in the case $d=0$. A vertical slice of the surface gives a particular density-flux assignment for a $b$ value.}\label{fig:flux}
\end{figure}

Now, substituting the measure $\Gamma^{\theta(v)}$ of \eqref{eq:statmarginal} into the previous display we get
\begin{equation*}
G(v) = \frac{b\cdot(1 - 2b)\cdot(\cosh(\theta(v)) + d \cdot \sinh(\theta(v))) + 2 b^2}{(1 + 2b\cdot(\cosh(\theta(v)) - 1))^2}.
\end{equation*}
Note that $G$ is \emp{even} (i.e.\ symmetric to the origin) if and only if $d=0$. This case is shown in Figure \ref{fig:flux}. In the upcoming sections we will only investigate the case when $d=0$.

Now, taking into consideration \eqref{eq:densityparam} we can realize the density parametrization of the flux. Elementary calculus in combination with hydrodynamic results as described in Section \ref{sec:hydro} then yield
%%%%%%%%%%%%%%%%%%%%%%%%%%%%%%%%%%%%%%%%
%%%%%%%%%%%%%%%%%%%%%%%%%%%%%%%%%%%%%%%%
\begin{theorem}\label{thm:main}
 The hydrodynamic limit as described in Section \ref{sec:hydro} applies to our model \eqref{eq:infgen} with rates \eqref{eq:modelrates} in the attractive range $c\leq\frac{1-|d|}2$ with hydrodynamic flux $H:[-1, 1]\to[0,+\infty)$ given by
\begin{equation*}
H(v) = \frac{v\cdot(d - v)}{2} +
\frac{(1 - d\cdot v)\cdot(4 b^2 + (1-2b)^2\cdot v^2)}{2 \,\big(4 b^2 + (1 - 2 b) \cdot\sqrt{4 b^2 + (1 - 4 b)\cdot v^2}\big)}
\qquad (v\in[-1,1]),
\end{equation*}
where $b = \frac{1}{2 + c^{-1/2}}$ and $c,d > 0$ (see \eqref{eq:paramb}). In the case $d=0$, we simply denote the flux by $G$:
\begin{equation}\label{eq:fluxsym}
G(v) = -\frac{v^2}{2} +
\frac{4 b^2 + (1-2b)^2\cdot v^2}{2\,\big(4 b^2 + (1 - 2 b)\cdot\sqrt{4 b^2 + (1 - 4 b)\cdot v^2}\big)} \qquad (v\in[-1,1]).
\end{equation}
For $G$ the following hold: if
$b \in \big(\frac{1}{6},\,\frac{1}{2+\sqrt{2}}\big]\Longleftrightarrow c \in \big(\frac{1}{16},\,\frac{1}{2}\big]$
($b = \frac{1}{6}\;\Longleftrightarrow\; c=\frac{1}{16}$), then it is strictly (non-strictly) concave in $[-1,1]$. Otherwise when
$b \in \big(0,\, \frac{1}{6}\big)\Longleftrightarrow c \in \big(0,\, \frac{1}{16}\big)$
it is neither concave nor convex having two inflection points
\begin{equation}\label{eq:inflpoints}
\pm v_{\infl} = \pm\sqrt{\frac{(2 b^2 (1-2b))^{2/3} - 4 b^2}{1-4b}},
\end{equation}
which separate a concave-convex ($-v_{\infl}$) and a convex-concave ($v_{\infl}$) region, respectively.
\end{theorem}

\section{Entropy solutions}\label{sec:entropysols}%

We find the entropy solutions of the Riemann problem \eqref{eq:hydrogeneral} with step initial datum \eqref{eq:stepinitcond} for all possible pairs of initial densities $\ul,\ur\in[-1,1]$. Recalling Section \ref{sec:hydro}, for purely concave flux two valid scenarios can only happen: the discontinuity at zero is preserved over time (\emp{shock wave}) or it immediately smears out (\emp{rarefaction fan}). As we have seen the flux $G$ can be non-concave in some parameter region and this results in plenty of qualitatively different solutions for \eqref{eq:hydrogeneral} to be discussed in more details below.

\begin{figure}[!htb]
\centering
\includegraphics[width=0.65\paperwidth]{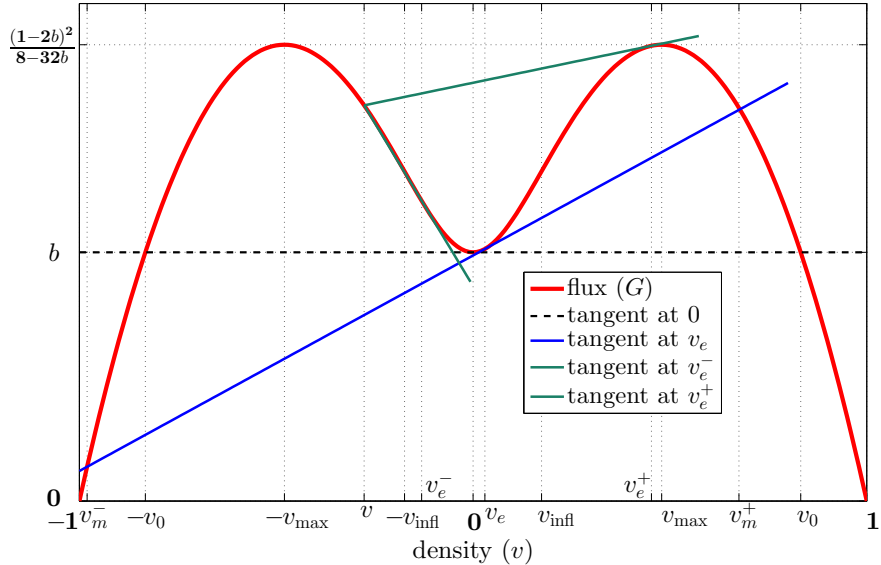}%
\caption[Flux (non-concave case)]{Hydrodynamic flux $G$ in the non-concave case and its tangents, where $b=0.08$ and $d=0$. From a general point $v\in[-1,1]$ one can possibly draw at most two lines that are tangential to $G$ (see $v_e^{-}$ and $v_e^{+}$). On the other hand a tangent line at $v_e$ can have two intersections with the graph of $G$ (see $v_m^{-}$ and $v_m^{+}$). $\pm v_{\max}$ are the maximum points while $\pm v_{0}$ are the intersections of $G$ with the constant (dashed) line at $b$.}\label{fig:tangent}
\end{figure}

\paragraph{Auxiliary calculations.} Based on Theorem \ref{thm:main}, we introduce some further notations being valid \emp{only} in the range $b\in(0,\frac{1}{6})$. The positive global maximum point of \eqref{eq:fluxsym} is given by
\[
v_{\max} = \frac{1}{2}\sqrt{\frac{(1 + 2b)(1 - 6b)}{1 - 4b}},
\]
while the negative one is at $-v_{\max}$. The maximum value is then $G(v_{\max}) = G(-v_{\max}) = (1 - 2b)^2 / (8 - 32b)$.
At $0$ the tangent line is horizontal and $G$ has a local minimum with $G(0) = b$. This tangent intersects the curve in two points symmetrically to the origin. The positive one is at
\begin{equation*}
v_0 = \sqrt{\frac{(1 - 2b)(1 - 6b)}{1 - 4b}},
\end{equation*}
while the negative one is at $-v_0$. See Figure \ref{fig:tangent}.
%%%%%%%%%%%%%%%%%%%%%%%%%%%%%%%%%%%%%%%%
%%%%%%%%%%%%%%%%%%%%%%%%%%%%%%%%%%%%%%%%

To determine the lower convex and upper concave envelopes of $G$ for different densities (recalling \eqref{eq:generalentropysol}) and to finally obtain the phase diagram we need to address some further special points of the flux $G$. Below we essentially define the tangential points which will be denoted by $v_e$'s whereas $v_m$'s will indicate the intersection points of $G$ with a tangent line. In more details the tangent line drawn from the point $(v_e, G(v_e))$ intersects the graph $G$ at (another) point $(v,G(v))$ if and only if it satisfies
\begin{equation}\label{eq:tangentpoint}
G(v) = G(v_{e}) + G'(v_{e})(v - v_{e}).
\end{equation}
We can look at \eqref{eq:tangentpoint} in two different ways.
\begin{itemize}[leftmargin=*]
\item One can look for an intersection $(v_m,G(v_m))$ of the graph and the tangent line touching at a \emp{given} $(v_{e}, G(v_e))$. The intersection that lies on the left (right) of $v_{e}$ is denoted by $v_m^{-}(v_{e})$ [$v_m^{+}(v_{e})$]. If $v_m^{-}(v_{e})$ [$v_m^{+}(v_{e})$] did not exist then $v_m^{-}(v_{e}):\,=-\infty$ [$v_m^{+}(v_{e}):\,=+\infty$]. See Figure \ref{fig:tangent}.
\item From a \emp{fixed} $(v,G(v))$ one can solve \eqref{eq:tangentpoint} for a tangent point $(v_{e},G(v_e))$ of $G$. If there exist more than one solutions then let $v_{e}^{-}(v)$ be the closer while let $v_{e}^{+}(v)$ be the farther one from $v$. By default $v_{e}(v):\,=v_{e}^{-}(v)$. If only one solution exits then $v_{e}(v)=v_{e}^{-}(v)=v_{e}^{+}(v)$. See Figure \ref{fig:tangent}.
\end{itemize}
In the following the previously defined points, illustrated by Figure \ref{fig:tangent}, are heavily used to determine those regions where the dynamics follows a shock wave or a rarefaction fan.
\paragraph{Phase diagram.} Notice that if $G''$ does not change sign in $[\min(\ul,\ur),\max(\ul,\ur)]$, furthermore the secant of $G$ connecting $(\ul,G(\ul))$ with $(\ur,G(\ur))$ is located \emp{below} (\emp{above}) its graph and
\begin{itemize}[leftmargin=*]
\item $\ul < \ur$ ($\ul > \ur$), then the characteristic lines `converge' and collide resulting in the formation of a \emp{shock wave}.
\item $\ul > \ur$ ($\ul < \ur$), then the characteristic lines `diverge' and the initial discontinuity smoothens out according to the physically relevant entropy solution, i.e.\ a \emp{rarefaction fan} forms.
\end{itemize}
Our first task is to explore these simple regions. Then in the general case we find the lower convex as well as upper concave envelopes of $G$ which will give shock and rarefaction parts of the solution via formula \eqref{eq:generalentropysol}. These functions can be expressed in a rather simple way by using the previously defined $v_e$'s and $v_m$'s.
%%%%%%%%%%%%%%%%%%%%%%%%%%%%%%%%%%%%%%%%
\begin{figure}[!htb]
\centering
\includegraphics[width=0.7\paperwidth]{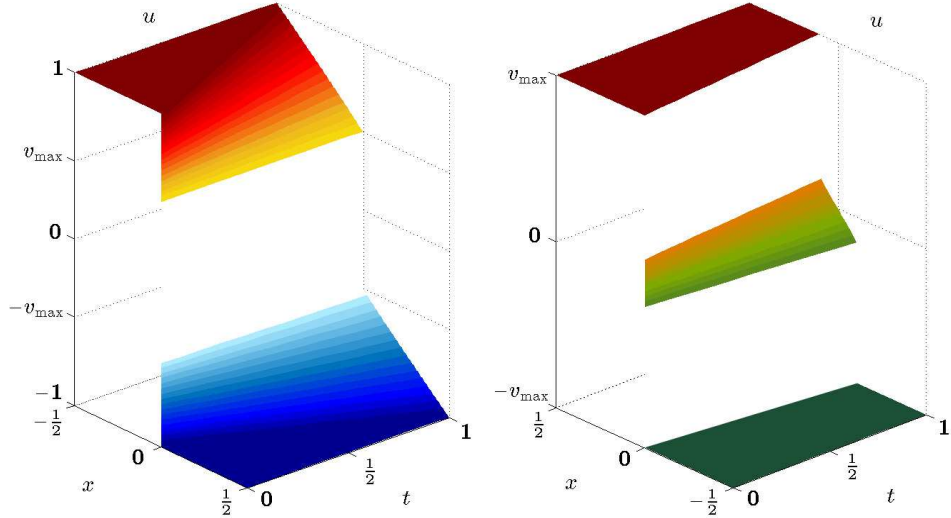}%
\caption[RSR and SRS]{Rarefaction fan -- Shock wave -- Rarefaction fan (left) and Shock wave -- Rarefaction fan -- Shock wave (right) profiles evolving in space-time when $\ul=1=-\ur$ and $\ul=-v_{\max}=-\ur$, respectively, where $b=0.08$ and $d=0$.}
\label{fig:rsr-srs}
\end{figure}
%%%%%%%%%%%%%%%%%%%%%%%%%%%%%%%%%%%%%%%%

The general case is as follows. Observe that if $b\geq\frac{1}{6}$ then shock wave (rarefaction fan) forms when $\ul < \ur$ ($\ul > \ur$). Henceforth we can assume that $b$ belongs to $(0,\frac{1}{6})$. For the simpler parts of the phase diagram we can apply the above rule. Recall the inflection points of \eqref{eq:inflpoints}.
\begin{itemize}[leftmargin=*]
\item[] \textbf{Shock wave (\textrm{S}) regions}. If $\ul < \ur$, then $-1\leq \ul \leq v_{m}^{-}(v_{e}(1))$ and $\ur \geq v_{m}^{+}(v_{e}(\ul))$; $v_{m}^{+}(v_{e}(-1)) \leq \ur \leq 1$ and $\ul \leq v_{m}^{-}(v_{e}(\ur))$; $-1\leq \ul\leq -v_{\infl}$ and $\ur\leq v_e(\ul)$; $v_{\infl}\leq \ur \leq 1$ and $\ul \geq v_e(\ur)$. If $\ul > \ur$, then $\ul,\ur\in[-v_{\infl},v_{\infl}]$; furthermore $v_{\infl} \leq \ul \leq v_{\max}$ ($-v_{\max}\leq \ur\leq -v_{\infl}$) and $v^{+}_e (\ul) \leq \ur \leq v^{-}_m(\ul)$ ($v^{+}_m(\ur)\leq \ul \leq v^{+}_e(\ur)$).
\item[] \textbf{Rarefaction fan (\textrm{R}) regions}. If $\ul < \ur$, then $\ul,\ur\in[-v_{\infl},v_{\infl}]$. On the other hand, if $\ul > \ur$, then $\ul,\ur\in[-1,-v_{\infl}]$ or $\ul,\ur\in[v_{\infl},1]$.
\end{itemize}

Now, it becomes a much more complex situation if the secant crosses the graph of $G$ above the interval $(\min(\ul,\ur),\max(\ul,\ur))$. We will discuss these in the following.
\begin{itemize}[leftmargin=*]
\item[] \textbf{Rarefaction fan -- Shock wave (\textrm{RS}) regions}. If $\ul < \ur$, then $v_{\infl} < \ur \leq 1$ and $-v_{\infl} \leq \ul < v_{e}(\ur)$. On the other hand if $\ul > \ur$, then $-v_{\max} \leq \ur < v_{\infl}$ and $v_{e}^{+}(\ur) < \ul \leq 1$.
\item[] \textbf{Shock wave -- Rarefaction fan (\textrm{SR}) regions}. If $\ul < \ur$, then $-1 \leq \ul < -v_{\infl}$ and $v_{e}(\ul) < \ur \leq v_{\infl}$. On the other hand if $\ul > \ur$, then $-v_{\infl} < \ul \leq v_{\max}$ and $-1 \leq \ur < v_{e}^{+}(\ul)$.
\item[] \textbf{Rarefaction fan -- Shock wave -- Rarefaction fan (\textrm{RSR}) region}: $\ul > v_{\max}$ and $\ur < -v_{\max}$. The shock wave in the middle has zero speed and $\lim_{x\to0^+}u(x,t)-\lim_{x\to0^-}u(x,t) = -2v_{\max}$ for all $t>0$. See Figure \ref{fig:rsr-srs}.
\item[] \textbf{Shock wave -- Rarefaction fan -- Shock wave (\textrm{SRS}) region}: $-1 \leq \ul < -v_{\infl}$ and $v_{\infl} < \ur < v_{m}^{+}(v_{e}(\ul))$. %or equivalently $v_{\infl} < \ur \leq 1$ and $v_{m}^{-}(v_{e}(\ur)) < \ul < -v_{\infl}$.
    The shocks are connected by a rarefaction fan and they are moving away from each other. See Figure \ref{fig:rsr-srs}.
\end{itemize}
%%%%%%%%%%%%%%%%%%%%%%%%%%%%%%%%%%%%%%%%
\begin{figure}[!htb]
\centering
\includegraphics[width=0.6\paperwidth]{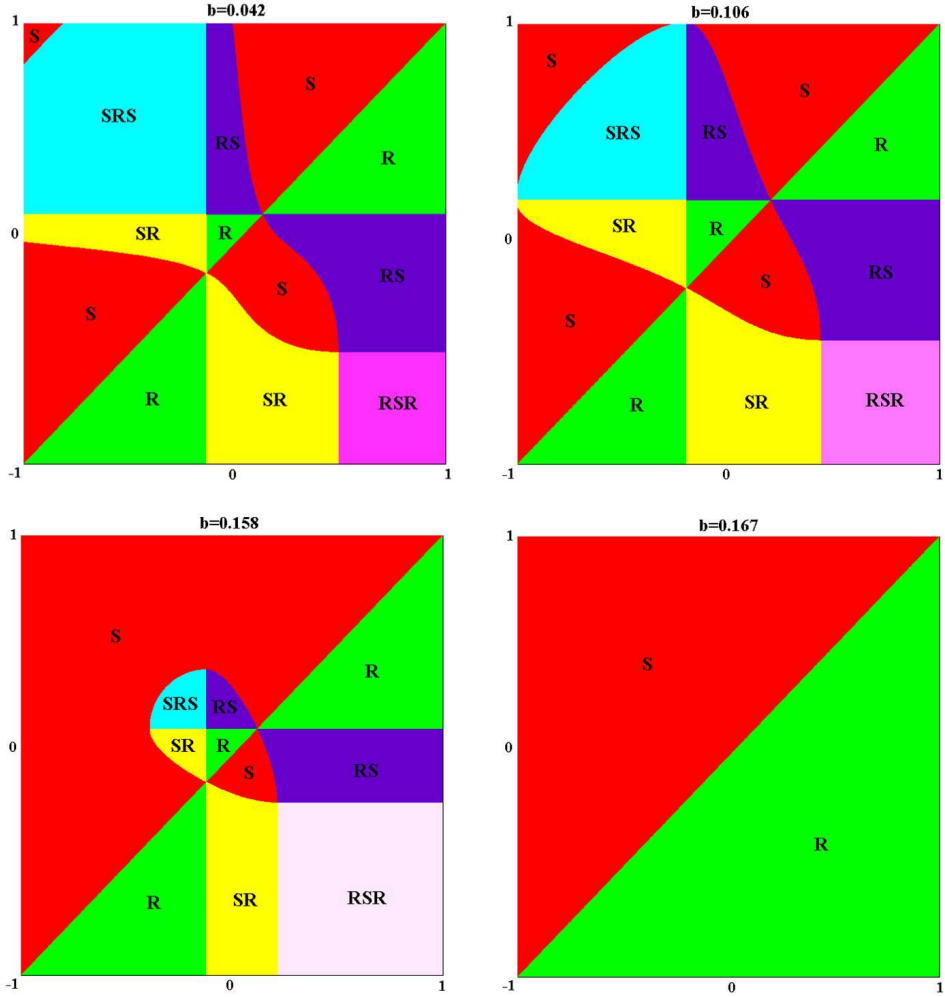}%
\caption[Phase diagram]{Phase diagram of the model with varying values of $b$, where $d=0$. The $x$ and $y$ axes represent the initial density on the left ($\ul$) and the right ($\ur$), respectively. $\mathbf{S}$ ($\mathbf{R}$) denotes a shock wave (rarefaction fan).}\label{fig:phasediagram}
\end{figure}
%%%%%%%%%%%%%%%%%%%%%%%%%%%%%%%%%%%%%%%%

The whole phase diagram can be seen in Figure \ref{fig:phasediagram} for four different $b$ values. The horizontal and vertical axes correspond to the initial density on the left ($\ul$) and the right ($\ur$), respectively. $\mathbf{S}$ denotes the formation of a shock wave, while $\mathbf{R}$ corresponds to a rarefaction fan. The model in the assumed $d=0$ case has a symmetry of swapping particles and antiparticles and flipping direction of the jumps at the same time. This results in the transformation $(\ul,\ur)\to (-\ur, -\ul)$ which shows as a symmetry on the diagrams (including reversing the order of $\mathbf{R}$'s and $\mathbf{S}$'s due to reflection of space). The other diagonal, $(\ul,\ur)\to (\ur,\ul)$ is \emp{not} a symmetry of the particle system due to positive particles jumping to the right, negatives to the left which fundamentally influences the hydrodynamic behavior.


\begin{thebibliography}{99}
\bibitem{bagurasa}
C.~Bahadoran, H.~Guiol, K.~Ravishankar and E.~Saada.
\newblock Euler hydrodynamics of one-dimensional attractive particle systems.
\newblock {\em Ann. Probab.}, 34(4):1339--1369, 2006.
%%
\bibitem{bagurasastrong}
C.~Bahadoran, H.~Guiol, K.~Ravishankar and E.~Saada.
\newblock Strong hydrodynamic limit for attractive particle systems on $\mathbb{Z}$.
\newblock {\em Electron. J. Probab.}, 15(1):1--43, 2010.
%%
\bibitem{pushasep}
A.~Borodin and P.~L.~Ferrari.
\newblock Large time asymptotics of growth models on space-like paths I: PushASEP.
\newblock {\em Electron. J. Probab.}, 13(50): 1380--1418, 2008.
%%
\bibitem{chowd}
D.~Chowdhury and J-S.~Wang.
\newblock Flow properties of driven-diffusive lattice gases: theory and computer simulation.
\newblock {\em Pyhs. Rev. E}, 65(4), 046126, 2002.
%%
\bibitem{coco}
C.~Cocozza-Thivent.
\newblock Processus des misanthropes.
\newblock {\em Z. Wahrsch. Verw. Gebiete}, 70(4):509--523, 1985.
%%
\bibitem{fossati}
M.~Fossati and L.~Quartapelle.
\newblock The Riemann Problem for Hyperbolic Equations under a Nonconvex Flux with Two Inflection Points.
\newblock 104 pages, 18{\`eme} Seminaire on la Mecanique des Fluides Numerique, Institut Henri Poincaré, France, January 2006. Public link: \url{http://arxiv.org/abs/1402.5906}.
%%
\bibitem{hager}
J.~S.~Hager, J.~Krug, V.~Popkov and G.~M.~Sh{\"u}tz.
\newblock Minimal current phase and universal boundary layers in driven diffusive systems.
\newblock {\em Phys. Rev. E}, 63(5), 056110, 2001.
%%
\bibitem{hayes}
B.~T.~Hayes and P.~G.~LeFloch.
\newblock Non-classical shocks and kinetic relations: scalar conservation laws.
\newblock {\em Arch. Rational Mech. Anal.}, 139(1):1--56, 1997.
%%
\bibitem{holdenrisebro}
H.~Holden and N.~H.~Risebro.
\newblock {\em Front Tracking for Hyperbolic Conservation Laws}.
\newblock Springer, 2011.
%%
\bibitem{kls}
S.~Katz, J.~L.~Lebowitz and H.~Spohn.
\newblock Nonequilibrium steady states of stochastic lattice gas models of fast ionic conductors.
\newblock {\em J. Stat. Phys.}, 34(3):497--537, 1984.
%%
\bibitem{kipnislandim}
C.~Kipnis and C.~Landim.
\newblock {\em Scaling Limits of Interacting Particle Systems}.
\newblock Springer, 1999.
%%
\bibitem{landim}
C.~Landim.
\newblock Conservation of local equilibrium for attractive particle systems on $\mathbb{Z}^d$.
\newblock {\em Ann. Probab.}, 21(4):1782--1808, 1993.
%%
\bibitem{liggettbible}
T.~M.~Liggett.
\newblock {\em Interacting Particle Systems}.
\newblock Springer, 1985.
%%
\bibitem{popkovetal}
V.~Popkov, A. R{\'a}kos, R.~D.~Willmann, A.~B.~Kolomeisky and G.~M.~Sch{\"u}tz.
\newblock Localization of shocks in driven diffusive systems without particle number conservation.
\newblock Phys. Rev. E 67(6), 066117, 2003.
%%
\bibitem{popkov}
V.~Popkov and G.~M.~Sch{\"u}tz.
\newblock Steady-state selection in driven diffusive systems with open boundaries.
\newblock {\em Europhys. Lett.}, 48(3), 1999.
%%
\bibitem{rezakhanlou}
F.~Rezakhanlou.
\newblock Hydrodynamic limit for attractive particle systems on $\mathbb{Z}^d$.
\newblock {\em Comm. Math. Phys.}, 140(3):417--448, 1991.
%%
\bibitem{smoller}
J.~Smoller.
\newblock {\em Shock Waves and Reaction-Diffusion Equations}.
\newblock Springer, 1983.
%%
\bibitem{junction}
B.~Tian, R.~Jiang, Z.-J.~Ding, M.-B.~Hu and Q.-S.~Wu.
\newblock Phase diagrams of the Katz--Lebowitz--Spohn process on lattices with a junction.
\newblock {\em Phys. Rev. E}, 87(6), 062124, 2013.
%%
\bibitem{tothvalkoperturb}
B.~T\'oth and B.~Valk\'o.
\newblock Between equilibrium fluctuations and Eulerian scaling: perturbation of equilibrium for a class of deposition models.
\newblock {\em J. Stat. Phys.}, 109(1):177--205, 2002.
%%
\bibitem{riemann}
I.~T\'oth.
\newblock Hydrodynamic behavior of an interacting particle system (in Hungarian). Master's Thesis,
\newblock Budapest University of Technology and Economics, Department of Stochastics, 2005.
%%
\bibitem{sim}
I.~T\'oth.
\newblock Simulation of a totally asymmetric attractive interacting particle system.
\newblock Public link: \url{http://github.com/tothi/riemann-nonconvex/}, 2016.
\end{thebibliography}
\end{document}